\documentclass[12pt]{article}
\usepackage{latexsym}
\usepackage{amsfonts}
\usepackage{amsthm}
\usepackage{amsmath}
\usepackage[dvips]{graphics}   
\begin{document}
\title{Groups which are not  properly 3-realizable}
\author{Louis Funar$^1$, Francisco F. Lasheras$^2$  and  Du\v{s}an  Repov\v{s}$^3$
\footnote{
{{Emails}: funar@fourier.ujf-grenoble.fr (L.Funar), 
lasheras@us.es (F.F.Lasheras), dusan.repovs@guest.arnes.si (D.Repov\v{s})}}
 \\
{\scriptsize {$^1$\em  Institut Fourier BP 74, UFR Math\'ematiques, Univ.Grenoble I 38402 Saint-Martin-d'H\`eres Cedex, France}} \\
{\scriptsize {$^2$\em Departamento de Geometria y Topologia, Universidad de Sevilla,  
Apdo 1160,  41080 Sevilla, Spain}}\\
{\scriptsize {$^3$\em Faculty of Mathematics and Physics, University 
of Ljubljana, P.O. Box 2964, Ljubljana 1001, Slovenia}}
}
\date{\today}
\maketitle
{\abstract
A group is properly 3-realizable if it is 
the fundamental group of a  compact polyhedron whose universal covering 
is proper homotopically equivalent to some 3-manifold. 
We prove that when such a group is also quasi-simply filtered then it   
has {\em pro-(finitely generated free) fundamental group at infinity} 
and {\em semi-stable ends}.  
Conjecturally the quasi-simply filtration assumption is superfluous. 
Using these restrictions we provide the first 
examples of finitely presented groups which are not properly 
3-realizable, for instance large families of Coxeter groups. 
}

\newtheorem{definition}{Definition}[section]

\newtheorem{lemma}{Lemma}[section]

\newtheorem{proposition}{Proposition}[section]

\newtheorem{theorem}{Theorem}[section]

\newtheorem{remark}{Remark}[section]

\newtheorem{corollary}{Corollary}[section]

\newtheorem{ex}{Example}[section]

\newtheorem{conjecture}{Conjecture}

\newcommand{\Z}{{\mathbb{Z}}}      
\newcommand{\Q}{{\mathbb{Q}}}      
\newcommand{\R}{{\mathbb{R}}}     

\vspace{0.3cm}
{\em AMS Math. Subj. Classification}(2000): 57 M 50, 57 M 10, 57 M 30. 

\vspace{0.3cm}

\section{Introduction}
The aim of this paper is  to obtain  necessary conditions 
for a finitely presented group to be properly 3-realizable, which lead 
conjecturally to a  complete characterization.  
Lasheras introduced and studied this class of groups in 
\cite{CLR,CL,La}. Recall that:

\begin{definition}
A finitely presented group $\Gamma$ is said to be properly 3-realizable 
(abbreviated {\em P3R} from now on) if there exists 
a  compact 2-dimensional polyhedron $X$ with fundamental group 
$\Gamma$  such that the universal covering $\widetilde{X}$  
is proper homotopy equivalent to a 3-manifold $W^3$. 
\end{definition}

\vspace{0.1cm}\noindent 
Hereafter we  will consider only {\em infinite} groups $\Gamma$ and 
thus the associated 3-manifolds $W^3$ appearing in the definition above 
will be {\em non-compact}.  Notice that, in general,  
the 3-manifolds $W^3$  will also have {\em non-compact} boundary.

\begin{remark}
In the definition of a  P3R group one does not claim that the universal 
covering of  {\em any}  
compact 2-dimensional polyhedron $X$ with fundamental group $\Gamma$ 
is proper homotopy equivalent to a 3-manifold.  
However it was proved in (\cite{ACLQ},  Proposition 1.3)   
that  given a P3R group $G$ 
then for {\em any} 2-dimensional  compact polyhedron 
$X$ with fundamental group $G$,  
the universal covering of the wedge $X\bigvee S^2$ 
is proper homotopy equivalent to a 3-manifold. 
\end{remark}

\vspace{0.2cm}
\noindent 
Recall the following classical theorem about embeddings up to homotopy, 
due to Stallings. 
Let $P$ be a finite CW-complex of dimension $k$, let $M$ be a PL-manifold of 
dimension $m$ and let $f\colon P\to M$ be a $c$-connected map. 
If $m-k\geq 3$ and if $c\geq 2k-m+1$ then there exist a compact 
subpolyhedron $j\colon Q\hookrightarrow M$ and a homotopy equivalence $h\colon P\to Q$ such that 
$jh$ is homotopic to $f$. 
This was generalized to the non-compact 
situation in \cite{CFLQ} by replacing the connectivity with the proper 
connectivity. Recall that a locally finite CW complex is said to be 
properly $c$-connected  ($c\geq 1$) 
if its proper homotopy type can be represented 
by a CW complex whose $c$-skeleton is reduced to an end-faithful tree 
(see \cite{CFLQ} for details).    
Thus the proper homotopy type of a locally finite CW-complex 
$X$ of dimension $n$ is represented by a closed subpolyhedron 
of ${R}^{2n-c}$ if $X$ is properly $c$-connected.

\vspace{0.2cm}
\noindent 
In particular,  the universal covering 
$\widetilde X$ of an {\em arbitrary} compact 2-polyhedron $X^2$ is 
proper homotopy equivalent to a 4-manifold, because 
any 2-polyhedron embeds, up to {\em proper homotopy}, 
into $\R^4$.  Therefore P3R groups are singled out among the set of all finitely presented groups 
by the fact that the universal covering  $\widetilde X$  
of some compact polyhedron $X$ with given 
$\pi_1(X)$  is proper homotopy equivalent to a particular 4-manifold, namely 
the product of a 3-manifold with an interval.

\begin{remark}
Fundamental groups of compact 3-manifolds 
are  obviously P3R, but there also exist  P3R groups which are {\em not} 
3-manifold groups.  For instance, any ascending 
HNN extension of a finitely presented group is P3R 
(\cite{La}, see  also other explicit examples in \cite{CL}). 
Moreover,  given any  infinite finitely presented groups $G$ and $H$, their 
direct product $G\times H$ is P3R (according to \cite{CLR}). Further 
amalgamated products of P3R groups (and HNN extensions) 
over finite groups yield P3R groups (see \cite{CLQR1}).
\end{remark}

\vspace{0.2cm}
\noindent  Let us introduce very briefly, for the sake of completeness, 
some end invariants of non-compact spaces which will 
be used in the sequel.  Standard references where these notions are 
studied in detail are 
\cite{BQ,MS}. 

\vspace{0.2cm}
\noindent  
Given the sequence of homomorphisms $A_{i-1}\leftarrow A_i$, called {\em bonding morphisms}, 
one builds the {\em tower of groups} $A_0\leftarrow A_1\leftarrow \cdots$.  
A {\em pro-isomorphism} between the towers 
$A_0\leftarrow A_1\leftarrow \cdots$ and 
$B_0\leftarrow B_1\leftarrow \cdots$ is given by two sequences of morphisms 
$ B_{j_{2n+1}}\to A_{i_{2n+1}}$ and $A_{i_{2n}}\to B_{j_{2n}}$ where 
$0=i_1 <j_1<j_2<i_2<i_3<j_3<j_4<i_4<\cdots$, which commute with the 
respective compositions of bonding morphisms in the two towers.  
A pro-isomorphism class of towers of groups is called a {\em pro-group}.

\begin{definition}
A pro-group is said to be {\em pro-(finitely generated free)}  if it has a representative 
tower in which all groups involved are finitely generated free groups. 
\end{definition}


\vspace{0.2cm}
\noindent  
It was shown in \cite{La} that if a pro-group is pro-(finitely generated
free) and has a representative tower with surjective bonding maps, then it has a representative
telescopic tower (i.e., a tower in which both conditions hold simultaneously).

\vspace{0.2cm}
\noindent  
Pro-groups arise in topology by means of towers associated to exhaustions of non-compact spaces. 

\begin{definition}
If $X$ is a polyhedron then a proper map $\omega:[0,\infty)\to X$ is called a {\em proper ray}. 
Two proper rays define the same {\em end} if their restrictions to the subset of natural numbers 
are properly homotopic. 

\vspace{0.1cm}\noindent 
An end is called {\em semi-stable} if every two proper rays defining 
this end are actually properly homotopic; one also says that the two rays define 
the same strong end. 

\vspace{0.1cm}\noindent 
A finitely presented group 
has {\em semi-stable ends} if there exists a compact polyhedron 
$X$ with  the given fundamental group whose universal covering has semi-stable ends.
\end{definition}

\vspace{0.2cm}
\noindent  
Given now a  proper base ray $\omega$ in $X$ and an exhaustion 
$C_1\subset C_2\subset \cdots \subset X=\cup_{i=1}^{\infty} C_i$ by compact 
subpolyhedra,   
we can associate a tower of groups 
\[ \pi_1(X,\omega(0))\leftarrow \pi_1(X-C_1,\omega(1))\leftarrow \cdots \]
where the bonding morphisms are induced, on the one hand, by the inclusions of spaces 
and on the other hand,  by  the change of  base points  which are slid along the ray $\omega$ 
restricted to integral intervals. 

\begin{definition}
The {\em (fundamental) pro-group at infinity} of $X$ 
based at $\omega$, denoted $\pi_1^{\infty}(X,\omega)$, is 
the pro-group associated to the tower of groups 
\[ \pi_1(X,\omega(0))\leftarrow \pi_1(X-C_1,\omega(1))\leftarrow \cdots \]
Two rays defining the same strong end yield isomorphic pro-groups. 
In particular, if the end is semi-stable, the pro-group at infinity is 
an invariant of the end, and  called the  {\em (fundamental) pro-group of the end}.
The end is called {\em simply connected at infinity (or $\pi_1$-trivial)} if the  associated 
pro-group is pro-isomorphic to a tower of trivial groups. 

\vspace{0.1cm}\noindent 
The {\em (fundamental) pro-group at infinity of a finitely presented group} 
is the pro-group at infinity of the universal covering of a compact 
polyhedron with the 
given fundamental group. This depends of course, on the base ray (and thus only on the end 
if it is semi-stable), but not on the the  particular compact polyhedron we chose. 
\end{definition}

\begin{remark}
There are alternative  equivalent definitions of the semi-stability, 
in particular  the 
one used in Siebenmann's thesis: an end is called  semi-stable 
if its fundamental pro-group has a representative 
tower with {\em surjective} bonding morphisms (see also \cite{HR}).  
For the sake of completeness 
we recall that an end is called {\em stable} if there exist some representative tower 
in which all bonding morphisms are {\em isomorphisms}. 
Examples of Davis (see \cite{DM}) show that the ends of universal coverings of 
finite complexes  might be  {\em not stable}, although it is not known  
whether they should be always semi-stable.  
Notice that sometimes in the literature one uses the 
terms $\pi_1$-stable, $\pi_1$-semi-stable etc. 
for the corresponding notions introduced above. 
As already observed above, we can infer from \cite{La} that a semi-stable end 
having pro-(finitely generated free) fundamental pro-group at infinity admits a representative
{\em telescopic tower} for that fundamental pro-group at infinity.
\end{remark}

\vspace{0.2cm}
\noindent
If a group has semi-stable ends then the universal covering of {\em any}  compact polyhedron 
with the given fundamental group has semi-stable ends. 
Although there exist spaces whose ends are not semi-stable, there are still 
no known examples 
of finitely presented groups  (i.e. universal coverings of compact polyhedra)  
without semi-stable ends (see also \cite{GM,Mi}).

\vspace{0.2cm}
\noindent 
The main source of examples of  P3R groups is the paper of  Lasheras 
(\cite{La}) where 
it is proved that  {\em a one-ended finitely presented group which is semi-stable and whose  
fundamental pro-group at infinity is pro-(finitely generated free) is P3R}.
In particular, any  one-ended finitely presented group $\Gamma$ which 
is simply connected at infinity (and hence automatically semi-stable at infinity)  
is P3R.  

\vspace{0.2cm}
\noindent 
We expect the following to be a complete characterization of this class of groups: 

\begin{conjecture}[3-dimensional homotopy covering conjecture]
A finitely presented group is P3R iff each one of its ends is 
semi-stable and has pro-(finitely generated free) fundamental pro-group.  
\end{conjecture}

\begin{remark}
In \cite{LaRo} the authors proved the sufficient part of the conjecture, 
namley that a finitely presented group whose ends are 
semi-stable and have pro-(finitely generated free) fundamental pro-groups 
is P3R. 
\end{remark}

\vspace{0.2cm}
\noindent 
In  this paper we give evidence in the favor of this 
conjecture,  by proving it in the case when 
the group under consideration satisfies an additional hypothesis related 
to the geometric simple connectivity. 
In order to explain this  we have to introduce, following Brick - Mihalik 
(\cite{BM}) and Stallings (\cite{St2}), 
the following tameness condition for groups and spaces.  

\begin{definition}
A space $X$ is  called {\em quasi-simply filtered } (abbreviated {\em qsf}) if for any compact $C\subset X$ there exists 
a connected and simply connected compact $K$ together with a 
map $f:K\to X$ such that $f(K)\supset C$ and $f|_{f^{-1}(C)}:f^{-1}(C)\to C$ 
is a homeomorphism. 

\vspace{0.2cm}
\noindent  
A finitely presented group $\Gamma$ is called qsf if there exists a  
(equivalently, for every) 
compact polyhedron $P$  with fundamental group $\Gamma$ 
such that the universal covering $\tilde{P}$ is qsf. 
\end{definition}

\vspace{0.2cm}
\noindent  
The condition qsf is a rather mild assumption on finitely presented groups. 
There are still no known examples of groups which do not have  the qsf property 
and most classes of known groups, as hyperbolic, semi-hyperbolic, 
automatic, tame combable etc., are qsf (see \cite{FO,MT2}).

\vspace{0.2cm}
\noindent 
We can now state  our main result: 

\begin{theorem}\label{group}
If a finitely presented group is P3R and qsf then all of its ends 
are semi-stable and have pro-(finitely generated free) 
fundamental group at infinity. 
\end{theorem}

\begin{remark}
We do not know whether all finitely presented groups which have semi-stable ends 
and pro-(finitely generated free) fundamental groups at each 
end are actually qsf. Notice that by a theorem of Wright (see \cite{Ge}, Theorem 16.5.6),  one-ended 
groups with stable end having an element 
of infinite order must be either simply connected at infinity or pro-$\Z$ 
at infinity. 
Thus they are P3R by the result of Lasheras cited  above.     
\end{remark}

\begin{remark}\begin{enumerate}
\item The homotopy covering conjecture implies the well-known covering conjecture 
in dimension 3 which states  that the universal 
covering of an irreducible closed 3-manifold $M^3$  with 
infinite fundamental group is  simply connected at infinity.  
In fact, the  universal covering $\widetilde M$ is an open contractible 3-manifold (thus one-ended)  
which is semi-stable and has pro-(finitely generated free) fundamental pro-group at infinity. 
This implies that there exists an exhaustion by compact submanifolds $C_i$ such that 
$\pi_1(\widetilde M-C_i)$ are finitely generated and free. 
Tucker's criterion from \cite{Tu2}  
implies that the manifold $\widetilde M$ is a missing boundary manifold and thus 
it is homeomorphic to ${\rm int}(N^3)$, for a suitable  compact 
3-manifold $N^3$ with boundary. By the contractibility of the universal 
covering, each component 
of $\partial N^3$ is homeomorphic to a 2-sphere and this implies that ${\rm int}(N^3)$ 
(and hence $\widetilde M$) is simply connected at infinity. 

\item 
Conversely, it is obvious that the universal covering conjecture implies the homotopy covering 
conjecture {\em for closed 3-manifold groups} because open 3-manifolds which 
are simply connected at infinity are semi-stable and have pro-(finitely generated free)  
pro-group at infinity (in fact a trivial pro-group!). 

\item 
Notice that the universal covering $\widetilde{X}$ of a compact 2-polyhedron $X$ can 
never be proper homotopy equivalent to an open (simply connected) 
3-manifold $M^3$.  
In fact, the Poincar\'e duality would give us that the third cohomology group
with compact support $H^3_c (\widetilde{X}) $ is isomorphic to  
$H^3_c (M) = H_0 (M) = \Z$, which is impossible, as $\dim(\widetilde{X})=2$.
\end{enumerate}
\end{remark}

\begin{remark}
Let us consider the universal covering  $\widetilde M^3$,  
of a 3-manifold $M^3$ with boundary.  If the boundary is a union of spheres 
then $\widetilde M$ is obtained from the universal covering of a closed 
3-manifold (obtained by capping off boundary spheres by balls)  by deleting  
a collection of disjoint balls.  Assume that the boundary is non-trivial i.e.  
not a union of 2-spheres.  Then $M^3$ is Haken and thus, by Thurston's theorem,  
it is a geometric 3-manifold. Let us moreover assume that $M^3$ is 
{\em atoroidal}, i.e. there are no $\Z\oplus \Z$ embedded in $\pi_1(M)$ other than peripheral 
subgroups coming from the boundary torus components.   
Then Thurston's geometrization 
theorem tells us that $M^3$ is hyperbolic. Therefore  
the universal covering $\widetilde M$  is obtained  
geometrically by deleting  a collection of horoballs from the hyperbolic 3-space.  
In particular the pro-group at infinity of $\widetilde M$ is pro-(finitely generated free)  and 
its ends are semi-stable. Thus the conjecture holds for fundamental groups 
of atoroidal 3-manifolds with non-trivial boundary.  A similar but more involved 
discussion shows that it  also holds  for all 3-manifolds with non-trivial boundary 
(since they are geometric). 
\end{remark}

\begin{remark}
The homotopy covering conjecture implies that all 
1-relator groups are P3R. This is already known for  1-relator 
finitely ended groups (see \cite{CLQR2}). 
In fact, 1-relator groups are semi-stable at infinity (see \cite{MT})  and it was proved in (\cite{CLQR2}, Proposition 2.7) 
that their pro-groups at infinity are pro-(finitely generated free). 
Notice that 1-relator groups are also qsf (see \cite{MT2}). 
Recently, Lasheras and Roy (\cite{LaRo}) have extended  the results of 
\cite{CLQR2} to a class of groups which contains all
1-relator groups. 

\vspace{0.2cm}
\noindent 
It is presently unknown (but quite plausible) that any finitely presented  
group which is  qsf, semi-stable and 
has pro-(finitely generated free) pro-groups at infinity  is 
P3R. 
\end{remark}

\vspace{0.2cm}
\noindent 
As an application  of Theorem \ref{group} we will obtain explicit examples 
of groups which are not $P3R$, as follows. 

\begin{theorem}\label{nonex}
Let $\Gamma$ be one of the following: 
\begin{enumerate}
\item the fundamental group of a finite non-positively 
curved complex which is a homology 
$n$-manifold ($n\geq 3$),  but not a topological manifold.
We further assume that the link of every vertex is a topological manifold.  
\item the right angled 
Coxeter group associated to a flag complex $L$ whose geometric 
realization is a closed combinatorial $n$-manifold ($n\geq 3$) and 
$\pi_1(L)$ is not a free group.   
\end{enumerate}
Then $\Gamma$ is not P3R. 
\end{theorem}

\vspace{0.2cm}
\noindent 
In particular many Coxeter groups are {\em not} P3R. 
Similar examples were announced by 
Cardenas.


\vspace{0.5cm} 
\noindent  
{\bf Acknowledgements.} The authors are indebted to Ross Geoghegan 
and Valentin Poenaru for useful discussions and comments and to an 
anonymous referee for simplifying the proofs.   
The first author was supported by 
the Proteus program (2005-2006), no 08677YJ and the ANR Repsurf: 
ANR-06-BLAN-0311. 
The second author was supported by the project MTM 2007-65726 and 
the third author was supported by the Proteus program (2005-2006), 
no 08677YJ.

\section{Proofs}

\subsection{Tameness criterion for non-compact 3-manifolds}
Recall that a polyhedron $P$ is called 
{\em weakly geometrically simply connected} ({\em wgsc}) 
if it admits an exhaustion 
by compact connected subpolyhedra $P_1\subset P_2\subset \cdots $
such that $\pi_1(P_n)=0$, for all $n$. 
The wgsc property for polyhedra is the piecewise-linear 
analogue of the geometric simple connectivity of open 
manifolds, namely the existence of a proper handlebody decomposition 
without index one handles.

\vspace{0.2cm}
\noindent 
It is proved in \cite{Fu2,FT} that an open 3-manifold proper homotopy 
equivalent to a  weakly geometrically simply connected polyhedron 
is simply connected at infinity. In this section we will extend 
this result to non-compact 3-manifolds.  

\vspace{0.2cm}
\noindent 
In the realm of manifolds with boundary the relevant tameness 
condition that will replace the simple connectivity at infinity is 
the following: 

\begin{definition}
A manifold $W$ is called a {\em missing boundary manifold} (also called 
{\em almost compact}) if there exists a compact manifold with boundary 
$M$ and a closed subset $A\subset \partial M$ 
of the boundary (not necessarily a subcomplex) 
such that $W$ is homeomorphic to $M-A$. 
\end{definition}

\vspace{0.2cm}
\noindent 
Interesting examples of manifolds which are not missing boundary 
manifolds can be 
found in \cite{ST,Tu1}.

\vspace{0.2cm}
\noindent 
We first introduce  a family of 3-manifolds  
which is, in some sense, the smallest one containing the 
missing boundary 3-manifolds and allowing manifolds to 
have infinitely many boundary components. 
These manifolds will be the proper analog of  
the open manifolds which are simply connected at 
infinity in the non-compact case.  

\vspace{0.2cm}
\noindent 
Before we proceed, let us recall that a compact $0$-dimensional subset $C$  
is said to be {\em tame} (or tamely embedded) in $\R^n$
if there exists a homeomorphism of $\R^n$ sending $C$  
into a subset of $\R\times \{0\}\subset \R^n$.  
It is well-known that perfect (i.e. without isolated points) 
compact $0$-dimensional separable topological 
spaces are homeomorphic to the Cantor space. 
Hence the tameness condition above is mostly relevant 
for Cantor subsets of $\R^n$. Notice that there exist wild Cantor sets 
in any $\R^n$, with $n\geq 3$, while Cantor sets in 
$\R^2$ are tame, by a classical theorem of Bing (\cite{Bing}). 

\begin{definition}
A  {\em standard model} 
is a 3-manifold with boundary $V$ constructed as 
follows. Let $\{B_i\}_{i\in I}$ be a collection of pairwise disjoint 3-balls 
in the interior ${\rm int}(B)$ of the 3-ball whose radii go to $0$ and 
whose limit set  $L$ is a tame $0$-dimensional subset 
disjoint from $\partial B$. 
Let $X\supset L$ be a tame $0$-dimensional subset of 
${\rm int}(B)$ which is disjoint from ${\rm int}(B_i)$, for all $i\in I$, and 
$T\subset \partial B\cup \cup_{i\in I} \partial B_i$. Then we put 
$V=B - (X\cup T\cup_{i\in I} {\rm int}(B_i))$.  
Manifolds of this form, where $T\cap \partial B=\emptyset$, were called 
{\em ragged cells} by Brin and Thickstun in (\cite{BT}, pp.9-10).  
\end{definition}

\vspace{0.1cm}\noindent
In order to simplify some arguments we will use in 
the sequel the fact that 
there are no fake homotopy disks in dimension 3, 
as the Poincar\'e conjecture has  been settled  by Perelman in 
\cite{Pe1,Pe2} 
(see  a detailed and self-contained exposition of  Perelman's proof 
in \cite{MoT}).

\begin{remark}\label{modele}
\begin{enumerate}
\item Open simply connected 3-manifolds $V$ which are simply connected 
at infinity can be described as the manifolds of the form  
$S^3- X$, where $X$ is a tame $0$-dimensional compact subset 
of $B^3$. Alternatively, $V$ can be written as an ascending union of 
compact simply connected submanifolds, i.e. disks - with - holes, 
by the Poincar\'e Conjecture (see \cite{FT,Wa}). 
\item 
A simply connected missing boundary 3-manifold $V$ is homeomorphic 
to $M-T$, where $M$ is a simply connected compact 3-manifold
and $T$ is a closed subset of $\partial M$ (see e.g. \cite{Wa}). 
By the Poincar\'e Conjecture 
there is a finite set of pairwise disjoint balls $B_i$, $i\in I$ 
such that $V=B - (\cup_{i\in I}{\rm int}(B_i)\cup T)$ and 
$T$ is a closed subset of $\partial B\cup \cup_{i\in I} \partial B_i$. 
Thus standard models $V$ with finite $I$ correspond precisely 
to simply connected missing boundary manifolds. 
Actually, any standard model can be obtained by making 
connected sums of (possibly infinitely many)  
simply connected missing boundary manifolds. 
\end{enumerate}
\end{remark}

\begin{remark}
\begin{enumerate}
\item Another characterization of standard models was given by Brin and Thickstun 
(see\cite{BT}, Full End Description Theorem (b), p.10), as follows.  
Modulo the Poincar\'e Conjecture, the set of simply connected end 1-movable 
3-manifolds coincides with that of standard models. 
In particular, 3-manifolds with semi-stable ends are homeomorphic to 
standard models.  
\item 
Cardenas announced as an application 
of the Brin-Thickstun structure theorem (\cite{BT}), 
that 1-ended groups which are 
P3R and semi-stable have actually 
pro-(finitely generated free) pro-group at infinity.  
\end{enumerate}
\end{remark}

\begin{remark}\label{gfdtgroup}
The boundary  of a standard model consists of 2-spheres and  open planar 
surfaces. Each end has pro-(finitely generated free) 
fundamental group at infinity. 
In fact, the complement of an unknotted ball in a 
1-ended standard model is homotopy equivalent to the complement  
of a finite graph, namely a holed handlebody. 
Thus its fundamental group is a finitely generated free group.   
Moreover, each end of a standard model is semi-stable.  
\end{remark}

\vspace{0.2cm}
\noindent 
The homotopy covering  conjecture admits an (\`a priori stronger) 
restatement  as follows:
 
\begin{conjecture}
Given a finitely presented P3R group, the  universal covering  
of some compact 2-dimensional polyhedron with this fundamental group is proper  
homotopy equivalent to a standard model. 
\end{conjecture}

\begin{remark}
The equivalence between the two conjectures stated in this paper 
is a consequence of the Brin-Thickstun structure theorem (\cite{BT}). 
Details are left to the reader. 
\end{remark}

\vspace{0.2cm}
\noindent 
The wgsc property is not so useful  anymore  if we consider the 
tameness of 3-manifolds {\em with boundary}.

\vspace{0.2cm}
\noindent 
The antecedent of the papers \cite{Fu2,FT} is the paper \cite{Po1} of Poenaru 
in which the geometric simple connectivity is already defined and used for  
non-compact manifolds with boundary.

\vspace{0.2cm}
\noindent 
Poenaru proved in \cite{Po1} that an open 3-manifold  
is simply connected at infinity if the product with a 
closed $n$-ball (for some $n\geq 2$) 
is a geometrically simply connected manifold with boundary. One might 
therefore expect 
that the analogous statement is true for the more general case  
of non-compact 3-manifolds. However, we will have to consider products 
of non-compact manifolds and disks, namely  manifolds with corners. 
It is thus natural to look for the piecewise-linear analogue of the 
geometric simple connectivity of manifolds with boundary. Specifically, we set: 

\begin{definition}\label{plgsc}
A polyhedron $P$ is said to be {\em pl-gsc} if it admits an exhaustion by compact 
connected subpolyhedra $P_1\subset P_2\subset \cdots$ 
such that $\pi_1(P_n)=0$ and 
$\pi_1(A, A\cap P_n)=0$, for every connected component $A$ of 
$\overline{P_{n+1}-P_n}$ and all $n$. 
Equivalently, the  map induced by inclusion 
$\pi_1(A\cap P_n)\to \pi_1(A)$ is surjective for all $A$, as above and all $n$. 
\end{definition} 

\vspace{0.2cm}
\noindent 
This definition is consistent with the previous ones since, by 
using Smale's theorem, a non-compact manifold of dimension 
$n\geq 6$ is pl-gsc iff it is gsc. 
Moreover the gsc and pl-gsc are equivalent for open manifolds 
without any dimensional restrictions.  
The pl-gsc is stronger than the wgsc 
for 3-manifolds with boundary.

\vspace{0.1cm}
\noindent
We start by recalling the following tameness condition for  
topological spaces, which is directly related to the qsf:

\begin{definition}
A non-compact PL space $X$ is called {\em Tucker} if the fundamental group 
of each component of $X-K$ is 
finitely generated, for any finite subcomplex $K\subset X$. 
\end{definition}

\vspace{0.1cm}
\noindent 
This definition was motivated by Tucker's work \cite{Tu2} on 3-manifolds, 
who proved that  
a $P^2$-irreducible connected 3-manifold is a  missing boundary 3-manifold
if and only if it is Tucker.

\vspace{0.2cm}\noindent 
The principal result of this section is the following extension of the 
result of \cite{Fu2,FT} to arbitrary non-compact 3-manifolds, as follows: 

\begin{proposition}\label{main}
A non-compact 3-manifold which has the proper homotopy type of a 
pl-gsc polyhedron is homeomorphic to a standard model. 
In particular, each end is 
semi-stable and its fundamental pro-group at infinity  
is pro-(finitely generated free). 
\end{proposition}
\begin{proof}
According to \cite{Fu2,FT} the interior ${\rm int} (W^3)$ 
is homeomorphic to a sphere minus a tame 0-dimensional subspace and 
in particular, it is simply connected at infinity. 

\vspace{0.2cm}\noindent 
Observe that a pl-gsc polyhedron is Tucker, by an easy application  
of the van Kampen theorem. Furthermore, if $W^3$ is proper homotopically equivalent to 
a Tucker polyhedron then $W^3$ is also Tucker (see \cite{O}).
 
\vspace{0.2cm}\noindent 
Let $e$ denote one of the countably many ends of ${\rm int}(W^3)$ 
which in $W^3$ has the boundary $\partial_eW^3$ associated to it.   
Consider a partial  (Freudenthal) end-point compactification of 
${\rm int}(W^3)$ which closes off all its ends but $e$. Recall that 
the end-point compactification of $X$ is a connected 
space $\widehat{X}$ containing $X$ as an open subset, 
with $\widehat{X}-X$ totally disconnected, such that 
for each $p\in \widehat{X}-X$, $U$ a connected open neighborhood 
of $p$ in $\widehat{X}$ the set $U-(\widehat{X}-X)$ is connected. 
In other words, we have one compactification point for 
each end of $X$.  

\vspace{0.2cm}\noindent
The end-point compactification is a manifold at some end 
if and only if the end is simply connected at infinity (see \cite{E,Wa}). 
We therefore obtain a simply connected 3-manifold $Z^3_e$ with boundary 
$\partial_eW^3$ whose interior ${\rm int}(Z^3_e)$ has only one end. 
Therefore $Z^3_e$ is irreducible. 
Observe that $Z^3_e$ has also the Tucker property. Using 
Tucker's criterion from \cite{Tu2} we deduce that $Z^3_e$ is a 
missing boundary manifold, and thus of the form 
$B^3 - T_e$,  where $T_e$ is a closed subset of 
$\partial B^3$ and $B^3$ is a 3-ball. 

\vspace{0.2cm}\noindent 
Use this method for each end $e$ of ${\rm int}(W^3)$ having 
a boundary associated to it. We can recover $W^3$ as the 
intersection of all $Z^3_e$ punctured along a tame Cantor subset 
corresponding to those ends having no boundary associated to them. 
Alternatively, $W^3$ is an infinite 
connected sum of all $Z^3_e$ punctured along a tame Cantor subset.  
Therefore $W^3$ is homeomorphic to a standard model.   
\end{proof}

\begin{remark}\label{ph}
We say that  the polyhedron $M$ is properly homotopically dominated by the 
polyhedron $X$ if there exists a PL map $f:M\to X$ whose mapping cylinder 
properly retracts on $M$. 
Then a manifold $W^3$ which is 
properly homotopically dominated by a Tucker polyhedron is also  
Tucker (see \cite{O}). 
This shows that 
Proposition \ref{main} extends to 3-manifolds which are  properly homotopically 
dominated by a pl-gsc polyhedron. 
\end{remark}

\subsection{Proof of Theorem \ref{group}}

We first prove: 
\begin{proposition}
If   the finitely presented group $G$ is P3R and qsf then there exists a 2-polyhedron $X$ 
with fundamental group $G$ such that $\widetilde{X}$ is  pl-gsc {\em and} 
proper homotopy equivalent to a 3-manifold $W^3$. 
\end{proposition}
\begin{proof}
Since $G$ is qsf it follows that for  {\em any} polyhedron $Y$ with fundamental group 
$G$,  the universal covering $\widetilde Y$ is qsf (see \cite{BM}). Take $Y$ to be a closed 
5-manifold with fundamental group $G$.  Then $\widetilde{Y}$ is an open 
5-manifold. It was proved in  (\cite{FO}, Proposition 3.2) that any  open simply connected  manifold 
 of dimension at least 5 which is qsf is actually gsc,  
as a consequence of general transversality results. 
It follows that $\widetilde{Y}$ is gsc.  We  triangulate $Y$ and get an equivariant 
triangulation of $\widetilde{Y}$.  Then the triangulated $\widetilde{Y}$ is a 
pl-gsc polyhedron.  The pl-gsc property is preserved when passing to 
the 2-skeleton. This means that the 2-skeleton $Z$ of the triangulation of $Y$ has the 
property that $\widetilde{Z}$ is pl-gsc. 
 
\vspace{0.2cm} 
\noindent  
It was proved in (\cite{ACLQ},  Proposition 1.3), as an application of Whitehead's theorem,  
that  given a P3R group $G$, for {\em any} 2-dimensional  compact polyhedron 
$X$ of fundamental group $G$, the universal covering of the wedge $X\bigvee S^2$ 
is proper homotopy equivalent to a 3-manifold. In particular, 
this holds when taking 
the 2-polyhedron $Z$ from above and thus $\widetilde{Z\bigvee S^2}$ is 
proper homotopy equivalent to a 3-manifold. 
Moreover, $\widetilde{Z\bigvee S^2}$ is made of one copy of $\widetilde{Z}$ 
with infinitely many $S^2$'s attached on it. In particular, 
if  $\widetilde{Z}$ is pl-gsc then 
it is immediate that $\widetilde{Z\bigvee S^2}$ is also pl-gsc. 
Therefore $X=Z\bigvee S^2$ has the required properties. 
\end{proof}

\vspace{0.2cm} 
\noindent  
{\em End of the proof of Theorem \ref{group}}. Let assume that we have a group $G$ which is 
both P3R and qsf. The previous proposition shows that there exists some 2-polyhedron $X$ 
such that $\widetilde{X}$ is pl-gsc and also proper homotopy equivalent to some 3-manifold $W^3$. 
Looking the other way around  we can apply Proposition \ref{main} to 
the 3-manifold $W^3$ (since it is proper homotopy equivalent 
to a pl-gsc polyhedron) and obtain that $W^3$ is homeomorphic to the standard model. 
In particular, $W^3$ has semi-stable ends and its pro-groups at infinity 
are pro-finitely generated free, as claimed. By the proper homotopy invariance of these end invariants 
$\widetilde{X}$ has the same properties. This proves Theorem \ref{group}. $\Box$

\begin{remark}
It follows by Remark \ref{ph} that 
Theorem \ref{group} holds  for the qsf groups $G$ for which there exists 
a  finite complex $X$ such that $\widetilde{X}$ properly homotopically 
dominates some 3-manifold. In particular, these groups are $P3R$.  
\end{remark}

\subsection{Proof of Theorem \ref{nonex}}
First, recall that groups acting properly cellularly and co-compactly 
on a CAT(0)-complex are wgsc and qsf (see \cite{FO,MT2}). Thus Coxeter groups 
and fundamental groups of  finite non-positively curved complexes are qsf. 

\vspace{0.2cm} 
\noindent  
Let us consider a finite non-positively curved complex $X$. We will 
use the criterion for the semi-stability given  in \cite{BMM}, 
which also provides a way to understand the pro-group at infinity. 
The link of a vertex in $X$ can be given a piecewise spherical metric.
Let $p$ be a point of the link of some vertex. The set of points of the link  
which are at distance at least $\frac{\pi}{2}$ from $p$ is called 
the punctured link. The punctured link deformation retracts 
onto the maximal subcomplex of the link  that it contains.  
The main theorem of \cite{BMM} states that if the links and the punctured 
links of $X$ are connected then $\widetilde X$ is has a semi-stable end. 

\vspace{0.2cm} 
\noindent  
If $X$ is a homology $n$-manifold both the links and the punctured links  
have the same $k$-homology as the $(n-1)$-sphere, for $k\leq n-2$.
In particular they are connected. On the other hand, there is at least 
one vertex $v$ of $x$ whose link is not simply connected, 
since the complex $X$ is not a topological manifold. The fundamental group 
of the link is then perfect and non-trivial and thus it cannot be a free 
group. The  complement of a punctured link within the link is the set of 
points of distance at most $\frac{\pi}{2}$ from the puncture. Since the 
metric structure of the link is $CAT(1)$ this complement is convex. 
Since we assumed the links to be topological manifolds each  
complement is a topological ball.

\vspace{0.2cm} 
\noindent  
In \cite{BMM} the Morse subdivision of $\widetilde X$ was defined  
as a geodesic subdivision induced by adding the critical points of the 
distance to a fixed base point. Let $\widetilde{X}_{>r}$ be the 
maximal subcomplex contained in the complement of the ball of radius $r$ 
in the Morse subdivision of $\widetilde{X}$. Since the distance is a 
Morse function  on a CAT(0)-complex  and the links are connected 
it is proved in \cite{BMM} that the inverse system 
\[ \pi_1(\widetilde{X}_{>0}) \leftarrow \pi_1(\widetilde{X}_{>1})
 \leftarrow \pi_1(\widetilde{X}_{>2})\leftarrow \cdots \]
has surjective bonding maps i.e. the end is semi-stable. 
Take the base point to be a lift of the 
vertex $v$. Then $\widetilde{X}_{>0}$ 
deformation retracts onto the link of $v$, 
and thus the first term of the inverse system is a non-free group.
Further $\widetilde{X}_{>r}$ deformation retracts (along geodesics)  
onto its boundary $\partial\widetilde{X}_{>r}$. 
On the other hand, one gets $\partial\widetilde{X}_{>r+1}$ from 
$\partial\widetilde{X}_{>r}$ by iterative use of the following procedure:
replace the complement of a punctured link embedded in  
$\partial\widetilde{X}_{>r}$ by the respective punctured link. 
Our hypothesis implies that the boundary of the complement of a punctured 
link is simply connected and hence, by induction and use of Van Kampen 
we find that $\pi_1(\partial\widetilde{X}_{>r})$ is a 
free product of fundamental groups of links. Moreover,   
the non-free group $\pi_1(\widetilde{X}_{>0})$ is a factor of the free 
product decomposition of every $\pi_1(\partial\widetilde{X}_{>r})$. 
The previous description shows also that the  
bonding maps correspond to forgetting a number of factors of the free product.  
Therefore $\widetilde{X}$ has a semi-stable end which is not pro-free. 
Since $\pi_1(X)$ is qsf it follows from the main theorem that it cannot be 
P3R. 

\vspace{0.2cm} 
\noindent 
The second part follows along the same lines. The topology at infinity 
of Coxeter groups was described in \cite{DM}. Recall that the  right angled 
Coxeter group $W_L$ associated to the flag complex $L$ is generated 
by the vertices of $L$ and the relations correspond to commutativity 
of adjacent vertices and the fact that these generators are of order two. 
Moreover $W_L$ acts on the Davis complex properly and cellularly.  
The Davis complex is a flag cubical complex and thus a CAT(0) complex. 
Thus $W_L$ is qsf (see also \cite{MT2}). 

\vspace{0.2cm} 
\noindent 
There is a natural filtration of the end defined by iterated neighborhoods 
of some vertex (see \cite{DM}). If $L$ is a closed connected combinatorial 
manifold then $W_L$ has one semi-stable end and the inverse sequence of 
fundamental groups is  as follows (see also (\cite{Ge}, Theorem 16.6.1)): 
\[ G \leftarrow G*G\leftarrow G*G*G \leftarrow \cdots \]
where $G=\pi_1(L)$ and each bonding map is  a projection annihilating 
the last factor. Thus if $L$ has dimension at least 2 and $G$ is not free 
then the fundamental group at infinity is not pro-free. 
The main theorem implies then that $W_L$ cannot be P3R. 
This settles Theorem \ref{nonex}.

\begin{remark}
We can infer from Remark \ref{gfdtgroup} that the higher homotopy 
groups at infinity $\pi_k^{\infty}(W)$ vanish for any standard model 
$W$ and $k\geq 3$. In particular,  
this furnishes another  practical tool for proving that 
a qsf finitely presented group $G$ is not P3R. Notice however, that 
this is a consequence of the fact that ends are semi-stable and 
pro-(finitely generated free).  
\end{remark}

\begin{small}
\bibliographystyle{plain}

\end{small}

\end{document}